\documentclass[11pt,leqno]{article}
\usepackage{amssymb}
\usepackage[mathscr]{eucal}
\usepackage{amsmath,amssymb,latexsym,theorem,bbm}
\usepackage{color}
\usepackage{appendix}

\newcommand{\NN}{\mathbb{N}}

\newcommand{\RR}{\mathbb{R}}

\newcommand{\cA}{{\mathcal A}}
\newcommand{\cB}{{\mathcal B}}

\newcommand{\dd}{\mathrm{d}}
\newcommand{\ee}{\mathrm{e}}

\newcommand{\EE}{\operatorname{\mathbb{E}}}
\newcommand{\PP}{\operatorname{\mathbb{P}}}

\newcommand{\cov}{\operatorname{Cov}}

\renewcommand{\leq}{\leqslant}
\renewcommand{\geq}{\geqslant}

\newcommand{\proofend}{\hfill\mbox{$\Box$}}

\numberwithin{equation}{section}

\theoremstyle{change} \theorembodyfont{\em}
\newtheorem{Lem}{Lemma.}[section]

\newtheorem{Pro}[Lem]{Proposition.}

\newtheorem{Def}[Lem]{Definition.}

\theorembodyfont{\rm}
\newtheorem{Rem}[Lem]{Remark.}

\begin{document}

\begin{center}
 {\bfseries\Large
 Examples of random fields that can be represented \\ as space-domain scaled stationary Ornstein-Uhlenbeck fields} \\[5mm]
 {\sc\large
  M\'aty\'as $\text{Barczy}$}
\end{center}

\vskip0.1cm

\noindent
   Faculty of Informatics, University of Debrecen,
   Pf.~400, H--4002 Debrecen, Hungary.

\noindent e--mail: barczy.matyas@inf.unideb.hu (M. Barczy).


\renewcommand{\thefootnote}{}
\vspace*{0.2cm}

\vspace*{-5mm}


\begin{abstract}
We give some examples of random fields that can be represented as space-domain scaled
 stationary Ornstein-Uhlenbeck fields defined on the plane.
Namely, we study a tied-down Wiener bridge, tied-down scaled Wiener bridges, a Kiefer process and so called \ $(F,G)$-Wiener bridges.
\end{abstract}

\section{Introduction}

In this note, we give some examples of random fields that can be represented as space-domain scaled
 stationary Ornstein-Uhlenbeck fields defined on the plane by specifying the space and domain transformations
 in question explicitly as well.
Before turning to fields, we recall a well-known one-dimensional example that a Wiener bridge can be represented
 as a space-time scaled stationary Ornstein-Uhlenbeck process.
Namely, if \ $(W(t))_{t\geq 0}$ \ is a standard Wiener process, then \ $S(t):=\ee^{-\frac{t}{2}} W(\ee^t)$, \ $t\in\RR$, \
 defines a strictly stationary centered Gauss process \ $S=(S(t))_{t\in\RR}$ \ defined on the real line with
  \ $\cov(S(s),S(t)) = \ee^{-\frac{\vert t-s\vert}{2}}$, \ $s,t\in\RR$,
 \ see, e.g., Doob \cite{Doo} or Shorack and Wellner \cite[Exercise 9, page 32]{ShoWel}.
The process \ $S$ \ is known as a stationary Ornstein-Uhlenbeck process defined on \ $\RR$.
\ Then a Wiener bridge \ $W(t) - t W(1)$, $t\in[0,1]$, \ from \ $0$ \ to \ $0$ \ over the time interval \ $[0,1]$
 \ generates the same law on \ $C([0,1])$ \ as the space-time scaled stationary Ornstein-Uhlenbeck process
 \[
   \begin{cases}
       \sqrt{t(1-t)}S\left(\ln \left(\frac{t}{1-t}\right)\right)
           & \text{if $t\in(0,1)$},\\
        0  & \text{if \ $t=0$ \ or \ $t=1$,}
     \end{cases}
 \]
 see, e.g., Shorack and Wellner \cite[Exercise 10, page 32]{ShoWel},
 where \ $C([0,1])$ \ denotes the space of continuous real-valued functions defined on \ $[0,1]$.

In Barczy and Kern \cite{BarKer3} we presented a class of Gauss-Markov processes which can be represented as space-time
 scaled stationary Ornstein-Uhlenbeck processes defined on the real line by giving examples as well, such as scaled Wiener bridges,
  Ornstein-Uhlenbeck type bridges, weighted Wiener bridges and so-called F-Wiener bridges.

In what follows, let \ $\NN$ \ and \ $\RR_+$ \ denote the set of positive integers and non-negative real numbers, respectively,
 let \ $\cB(\RR)$ \ be the set of Borel sets of \ $\RR$, \ and for \ $s,t\in\RR$,
 \ let \ $s\wedge t$ \ and \ $s\vee t$ \  denote \ $\min(s,t)$ \ and \ $\max(s,t)$, \ respectively.
For a subset \ $D\subseteq\RR^2$, \ $C(D)$ \ denotes the space of continuous real-valued functions on \ $D$.

\begin{Def}
A zero-mean Gauss field \ $\{W(s,t) : s,t\in\RR_+\}$ \ with continuous sample paths almost surely and with covariance function
 \begin{align*}
   \cov(W(s_1,t_1),W(s_2,t_2))
      = (s_1\wedge s_2) (t_1\wedge t_2),
      \qquad s_1,s_2,t_1,t_2\in\RR_+,
 \end{align*}
 is called a standard Wiener field (or a bivariate Wiener process).
\end{Def}

Here by a Gauss field, we mean that for all \ $n\in\NN$ \ and \ $(s_i,t_i)\in\RR_+^2$, $i=1,\ldots,n$, \
 the random variable \ $(W(s_1,t_1), \ldots,W(s_n,t_n))$ \ is \ $n$-dimensional normally distributed.
By the property that \ $W$ \ has continuous sample paths almost surely,
 we mean that \ $\PP(\{\omega\in\Omega : \RR_+^2\ni (s,t) \mapsto W(s,t)(\omega) \;\; \text{is continuous}\})=1$,
 \ where \ $(\Omega,\cA,\PP)$ \ denotes an underlying probability space on which the random variables in question are defined.

\begin{Def}\label{Def_OU_field}
Let \ $\alpha,\beta,\sigma>0$. \
A zero-mean Gauss field \ $\{X(s,t) : s,t\in\RR\}$ \ with continuous sample paths almost surely and with covariance function
 \begin{align*}
   \cov(X(s_1,t_1),X(s_2,t_2))
      = \frac{\sigma^2}{4\alpha\beta}
         \exp\Big\{-\alpha\vert s_2 - s_1\vert - \beta \vert t_2 - t_1\vert\Big\},
      \qquad s_1,s_2,t_1,t_2\in\RR,
 \end{align*}
 is called a stationary Ornstein-Uhlenbeck field with parameters \ $(\alpha,\beta,\sigma)$.
\end{Def}

Here stationarity means that for all \ $n\in\NN$ \ and \ $(s_i,t_i)\in\RR^2$, $i=1,\ldots,n$, \
 the distribution of the random vector
 \[
   (X(s+s_1,t+t_1),X(s+s_2,t+t_2),\ldots,X(s+s_n,t+t_n))
 \]
 does not depend on \ $(s,t)\in\RR^2$.

Next, we present a well-known result that a stationary Ornstein-Uhlenbeck field can be represented as a space-domain scaled
 standard Wiener field, see, e.g., Walsh \cite[page 271]{Wal} or Baran, Pap and van Zuijlen \cite[(3.1)]{BarPapZui}.
In particular, this shows the existence of a stationary Ornstein-Uhlenbeck field.
For completeness, we give a proof of this statement as well.

\begin{Pro}\label{prop_stac_OU_field}
If \ $\alpha,\beta,\sigma>0$ \ and \ $\{W(s,t) : s,t\in\RR_+\}$ \ is a standard Wiener field,
 then the random field
 \[
   Z(s,t):=\frac{\sigma}{2\sqrt{\alpha\beta}} \,\ee^{-\alpha s - \beta t}
           W(\ee^{2\alpha s},\ee^{2\beta t}),
          \qquad s,t\in\RR,
 \]
 is a stationary Ornstein-Uhlenbeck field with parameters \ $(\alpha,\beta,\sigma)$.
\end{Pro}

\noindent{\bf Proof.}
It can be readily seen that \ $Z$ \ is a zero-mean Gauss field with continuous sample paths almost surely.
Further, for all \ $s_1,s_2,t_1,t_2\in\RR$,
 \begin{align*}
 \cov(Z(s_1,t_1),Z(s_2,t_2))
   & = \frac{\sigma^2}{4\alpha\beta} \,
       \ee^{-\alpha s_1 - \beta t_1} \ee^{-\alpha s_2 - \beta t_2}
           \cov( W(\ee^{2\alpha s_1},\ee^{2\beta t_1}), W(\ee^{2\alpha s_2},\ee^{2\beta t_2} ))\\
   & = \frac{\sigma^2}{4\alpha\beta}
      \ee^{-\alpha (s_1 + s_2) - \beta (t_1+t_2)}
      (\ee^{2\alpha s_1}\wedge \ee^{2\alpha s_2})
      (\ee^{2\beta t_1}\wedge \ee^{2\beta t_2}) \\
   & = \frac{\sigma^2}{4\alpha\beta}
      \ee^{-\alpha (s_1 + s_2 - 2(s_1\wedge s_2))} \ee^{-\beta (t_1+t_2 - 2(t_1\wedge t_2))}\\
   & = \frac{\sigma^2}{4\alpha\beta}
      \ee^{-\alpha\vert s_1-s_2\vert - \beta\vert t_1-t_2\vert},
 \end{align*}
 where we used that \ $x+y - 2(x\wedge y) = \vert x-y\vert$, $x,y\in\RR$.
\proofend

We remark that despite the fact that the covariance function of a stationary Ornstein-Uhlenbeck field
 is of product type, a stationary Ornstein-Uhlenbeck field can not be represented as a product
 of two independent one dimensional Ornstein-Uhlenbeck processes (since such a product process is not a Gauss field).

In the present note, we show that a tied-down Wiener bridge, tied-down scaled Wiener bridges,
 a Kiefer process and so called \ $(F,G)$-Wiener bridges can be represented as space-domain scaled
 stationary Ornstein-Uhlenbeck fields defined on the plane.
In Remark \ref{Rem_counter} we point out that the form of the space-domain transformations used for representing a
 tied-down Wiener bridge as a space-domain scaled stationary Ornstein-Uhlenbeck field is not appropriate in case of a bivariate Wiener bridge.
Note that, due to Proposition \ref{prop_stac_OU_field}, a tied-down Wiener bridge, tied-down scaled Wiener bridges,
 a Kiefer process and \ $(F,G)$-Wiener bridges can be represented as a space-domain scaled Wiener field as well,
 however a Wiener field is not stationary.
The presented results may be used later on to calculate the distribution of the supremum location
 of the studied random fields on compact subsets.
For such an application in dimension one, see Barczy and Kern \cite[Section 5]{BarKer3}.

\section{Bivariate and tied-down Wiener bridge}\label{Section2}

Let \ $\{W(s,t) : s,t\in\RR_+\}$ \ be a standard Wiener field.
Deheuvels et al. \cite[formulas (2.15) and (2.16)]{DehPecYor} introduced
 a bivariate Wiener bridge
 \[
   B(s,t):= W(s,t) - stW(1,1),\qquad s,t\in[0,1],
 \]
 and a tied-down Wiener bridge
 \begin{align*}
  B_*(s,t)&:= B(s,t) - sB(1,t) - tB(s,1)\\
          &\;= W(s,t) - sW(1,t) - tW(s,1) + stW(1,1),
              \qquad s,t\in[0,1].
 \end{align*}
Then \ $\{B(s,t) : s,t\in\RR_+\}$ \ and \ $\{B_*(s,t) : s,t\in\RR_+\}$ \ are zero-mean Gauss fields
 with continuous sample paths almost surely and with covariance functions
 \[
   \cov(B(s_1,t_1),B(s_2,t_2))
      = (s_1\wedge s_2)(t_1\wedge t_2)
          - s_1s_2t_1t_2,
 \]
 and
 \[
   \cov(B_*(s_1,t_1),B_*(s_2,t_2))
      = (s_1\wedge s_2 - s_1s_2)(t_1\wedge t_2 - t_1t_2)
 \]
 for \ $s_i,t_i\in[0,1]$, $i=1,2$, \ respectively,
 see Deheuvels at el. \cite[formulas (2.17) and (2.18)]{DehPecYor}.
Note that the bivariate Wiener bridge \ $B$ \ is zero on the line segments between \ $(0,0)$ \ and \ $(0,1)$,
 \ and between \ $(0,0)$ \ and \ $(1,0)$, \ and at the point \ $(1,1)$ \ as well; \ while the tied-down Wiener bridge
 \ $B_*$ \ is zero on the border of a square with vertices \ $(0,0)$, $(0,1)$, $(1,0)$, \ and \ $(1,1)$.
\ Note also that the covariance function of \ $B_*$ \ is of product type (i.e., a product of functions depending
 on \ $(s_1,s_2)$ \ and \ $(t_1,t_2)$, \ respectively), while the covariance function of \ $B$ \ is not of that kind.

\begin{Pro}\label{prop_tied_dow_W-bridge}
Let \ $\{X(s,t) : s,t\in\RR\}$ \ be a stationary Ornstein-Uhlenbeck field with parameters \ $(\frac{1}{2},\frac{1}{2},1)$ \
 represented as in Proposition \ref{prop_stac_OU_field}.
Then the random field
 \begin{align*}
  U(s,t)
    :=\begin{cases}
        \sqrt{st(1-s)(1-t)}\, X\left(\ln\left(\frac{s}{1-s}\right), \ln\left(\frac{t}{1-t}\right) \right),
              & \text{if \ $(s,t)\in(0,1)^2$,}\\
            0 & \text{if \ $s\in\{0,1\}$ \ or \ $t\in\{0,1\}$,}
      \end{cases}
 \end{align*}
 generates the same law on \ $C([0,1]^2)$ \ as a tied-down Wiener bridge \ $B_*$.
\end{Pro}

\noindent{\bf Proof.}
First, we check that both fields \ $U$ \ and \ $B_*$ \ are zero-mean Gauss fields on \ $[0,1]^2$ \ with continuous sample paths almost surely.
The only property that does not follow immediately is that \ $U$ \ has continuous sample paths almost surely.
We need to prove that
 \[
   \PP\big( \{ \omega\in\Omega : [0,1]^2\ni(s,t)\mapsto U(s,t)(\omega) \;\;\text{is continuous}\}\big)=1.
 \]
Note that
 \begin{align*}
    &\{ \omega\in\Omega : [0,1]^2\ni(s,t)\mapsto U(s,t)(\omega) \;\;\text{is continuous} \}\\
    & = \{ \omega\in\Omega : [0,1]^2\ni(s,t)\mapsto U(s,t)(\omega) \;\;\text{is continuous in every point}\;\; (s_0,t_0)\in[0,1]^2\;\;\text{with}\;\; s_0\vee t_0 <1\} \\
    &\phantom{=\;} \cap \{ \omega\in\Omega : [0,1]^2\ni(s,t)\mapsto U(s,t)(\omega) \;\;\text{is continuous in every point}\;\; (s_0,t_0)\in[0,1]^2\;\;\text{with}\;\; s_0\vee t_0 =1\} \\
    & =: A_1\cap A_2,
 \end{align*}
 hence it is enough to prove that \ $\PP(A_1)=1$ \ and \ $\PP(A_2)=1$.
\ For all \ $s,t\in(0,1)$,
 \begin{align}\label{help2}
 \begin{split}
  U(s,t)&= \sqrt{st(1-s)(1-t)}\, X\left(\ln\left(\frac{s}{1-s}\right), \ln\left(\frac{t}{1-t}\right) \right) \\
        & = \sqrt{st(1-s)(1-t)}
            \exp\Big\{ -\frac{1}{2}\ln\left(\frac{s}{1-s}\right) - \frac{1}{2}\ln\left(\frac{t}{1-t}\right) \Big\}
            W\left(\ee^{\ln\left(\frac{s}{1-s}\right)},\ee^{\ln\left(\frac{t}{1-t}\right)}\right) \\
        & = \sqrt{st(1-s)(1-t)}
            \sqrt{\frac{1-s}{s}\cdot \frac{1-t}{t} }
           W\left(\frac{s}{1-s},\frac{t}{1-t}\right) \\
        & = (1-s)(1-t) W\left(\frac{s}{1-s},\frac{t}{1-t}\right).
 \end{split}
 \end{align}
Since the mapping \ $(s,t)\mapsto \left(\frac{s}{1-s},\frac{t}{1-t}\right)$ \ is a continuous homeomorphism of the
 set \ $\{(s,t)\in[0,1]^2: s\vee t <1\}$ \ onto \ $\RR_+^2$, \ by \eqref{help2}, we have
 \ $A_1 = \{ \omega\in\Omega : \RR_+^2\ni(u,v)\mapsto W(u,v)(\omega) \;\;\text{is continuous}\}$, \
 and hence, using that \ $W$ \ has continuous sample paths almost surely, we get \ $\PP(A_1)=1$.

Now we turn to prove that \ $\PP(A_2)=1$.
\ Recall that
 \begin{align}\label{help_LIL_Wiener_field1}
   \limsup_{r\to\infty}
       \sup_{(x,y)\in D_r^{(p)}} \frac{\vert W(x,y)\vert}
                       {\sqrt{4r\ln(\ln(r))}}
        = 1 \qquad \text{a.s.,} \qquad \forall\,p>\frac{1}{2},
  \end{align}
 where
 \begin{align*}
   & D_r^{(p)} := \Big\{ (x,y)\in\RR_+^2 : xy\leq r, \; 0\leq x\leq r^p,\; 0\leq y\leq r^p \Big\},
 \end{align*}
 see, e.g., Theorem 1.12.3 in Cs\"org\H o and R\'ev\'esz \cite{CsorRev}.
Using that a continuous function takes the limits of sequences to limits of sequences, we get
 \ $\{ \omega\in\Omega : \lim_{s\vee t \uparrow 1} U(s,t)(\omega) = 0 \}\subseteq A_2$.
\ By \eqref{help2}, for \ $s,t\in(0,1)$,
  \begin{align*}
  U(s,t)
    = \sqrt{\frac{4\ln\left(\ln\left(  \left( \frac{s}{1-s} +1\right)\left(\frac{t}{1-t}+1\right)\right)\right)}
           { \left( \frac{s}{1-s} +1\right)\left(\frac{t}{1-t}+1\right) }}
      \frac{W\left(\frac{s}{1-s},\frac{t}{1-t}\right)}
              {\sqrt{4 \left( \frac{s}{1-s} +1\right) \left(\frac{t}{1-t}+1\right)
                 \ln\left(\ln\left(  \left( \frac{s}{1-s} +1\right) \left(\frac{t}{1-t}+1\right)\right)\right)}},
 \end{align*}
 where
 \[
    \left(\frac{s}{1-s},\frac{t}{1-t}\right)
       \in D^{(1)}_{ \left( \frac{s}{1-s} +1\right) \left(\frac{t}{1-t}+1\right)}.
 \]
Here
 \[
   \left( \frac{s}{1-s} +1\right) \left(\frac{t}{1-t}+1\right)
     = \frac{1}{(1-s)(1-t)}
     = \frac{1}{(1-s\wedge t)(1-s\vee t)}
     \geq \frac{1}{1- s\vee t}
     \to\infty
     \qquad \text{as \ $s\vee t\uparrow 1$.}
 \]
Hence, using \eqref{help_LIL_Wiener_field1} and \ $\lim_{h\uparrow \infty} \frac{1}{h}\ln(\ln(h)) = 0$,
 \ we have \ $\lim_{s\vee t\uparrow 1}U(s,t)=0$ \ almost surely, yielding \ $\PP(A_2)=1$.

To conclude, it is sufficient to check that the covariance functions of \ $U$ \ and \ $B_*$ \ coincide.
If \ $0<s_1\leq s_2< 1$ \ and \ $0< t_1\leq t_2< 1$ \ (which can be supposed without loss of generality),
 then
 \begin{align*}
     & \cov(U(s_1,t_1), U(s_2,t_2))\\
     & = \sqrt{s_1t_1(1-s_1)(1-t_1)s_2t_2(1-s_2)(1-t_2)}\\
     &\phantom{=\;}
        \times\exp\Big\{-\frac{1}{2}\left( \ln\left(\frac{s_2}{1-s_2}\right) -  \ln\left(\frac{s_1}{1-s_1}\right)\right)
                 - \frac{1}{2}\left(\ln\left(\frac{t_2}{1-t_2}\right) -  \ln\left(\frac{t_1}{1-t_1}\right) \right) \Big\}\\
     & = \sqrt{s_1s_2t_1t_2(1-s_1)(1-t_1)(1-s_2)(1-t_2)}
         \sqrt{\frac{s_1}{1-s_1}\cdot \frac{1-s_2}{s_2}}
         \sqrt{\frac{t_1}{1-t_1}\cdot \frac{1-t_2}{t_2}}\\
     & = s_1t_1(1-s_2)(1-t_2)\\
     &= \cov(B_*(s_1,t_1),B_*(s_2,t_2)),
 \end{align*}
 where we used that the function \ $(0,1)\ni x\mapsto \ln\left(\frac{x}{1-x}\right)$ \ is strictly monotone increasing.
\proofend

In the next remark we present an alternative way for checking that \ $U$ \ defined in Proposition \ref{prop_tied_dow_W-bridge}
 has continuous sample paths almost surely.
We will use this approach in the proof of Proposition \ref{Prop_scaled_Wiener_field}.

\begin{Rem}\label{Rem_alternativ_proof}
Let \ $U$ \ be defined as in Proposition \ref{prop_tied_dow_W-bridge}.
Let \ $C_1:=\{ \omega\in\Omega : \RR_+^2\ni(s,t)\mapsto W(s,t)(\omega)\;\;\text{is continuous}\}$.
\ Then \ $\PP(C_1)=1$.
\ Recall that
\begin{align}
  & \limsup_{x\to\infty, \, y\to\infty}
           \frac{\vert W(x,y)\vert}
                {\sqrt{4xy\ln(\ln(xy))}}
        = 1 \qquad \text{a.s.}
  \label{help_LIL_Wiener_field2}
 \end{align}
 see, e.g., Theorem 1.12.2 in Cs\"org\H o and R\'ev\'esz \cite{CsorRev},
 and let \ $C_2$ \ be the set of those \ $\omega\in\Omega$ \ for which \eqref{help_LIL_Wiener_field2} holds
 (then \ $\PP(C_2)=1$).
\ If \ $(s,t)\to (s_0,t_0)\in[0,1]^2$ \ with \ $s_0\vee t_0<1$, \ then, by \eqref{help2},
 \ $U(s,t)(\omega)\to U(s_0,t_0)(\omega)$ \ for all \ $\omega\in C_1$.
\ If \ $(s,t)\to (s_0,t_0) = (1,1)$, \ then, by \eqref{help2} and \eqref{help_LIL_Wiener_field2}, we have
 \begin{align*}
 U(s,t)(\omega)
   = \sqrt{4st(1-s)(1-t) \ln\left(\ln\left(\frac{st}{(1-s)(1-t)}\right)\right)}
     \frac{W\left(\frac{s}{1-s}, \frac{t}{1-t} \right)(\omega)}
         {\sqrt{4 \frac{st}{(1-s)(1-t)} \ln\left(\ln\left(\frac{st}{(1-s)(1-t)}\right)\right)}}
     \to 0
 \end{align*}
 for all \ $\omega\in C_2$, \ where we used that \ $\lim_{h\uparrow \infty} \frac{1}{h}\ln(\ln(h)) = 0$.
\ If \ $(s,t)\to (1,t_0)$ \ with \ $t_0 \in [0,1)$, \ then, by \eqref{help2}, we have
 \begin{align*}
  U(s,t)
   & = s(1-t)\left(\frac{t_0}{1-t_0}+1\right) \sqrt{\frac{4\ln\left(\ln\left(\frac{s}{1-s}\left(\frac{t_0}{1-t_0}+1\right)\right)\right)}
           {\frac{s}{1-s}\left(\frac{t_0}{1-t_0}+1\right) }} \\
   &\phantom{=\;} \times \frac{W\left(\frac{s}{1-s},\frac{t}{1-t}\right)}
              {\sqrt{4\frac{s}{1-s}\left(\frac{t_0}{1-t_0}+1\right)\ln\left(\ln\left(\frac{s}{1-s}\left(\frac{t_0}{1-t_0}+1\right)\right)\right)}},
 \end{align*}
 where
 \[
    \left(\frac{s}{1-s},\frac{t}{1-t}\right)
       \in D^{(1)}_{\frac{s}{1-s} \left(\frac{t_0}{1-t_0}+1\right)}
 \]
 provided that \ $s$ \ is sufficiently close to \ $1$ \ (in fact, \ $s\in(1/2,1)$ \ is enough) and
 \ $t$ \ is sufficiently close to \ $t_0$.
\ By \eqref{help_LIL_Wiener_field1}, using again \ $\lim_{h\uparrow \infty} \frac{1}{h}\ln(\ln(h)) = 0$,
 \ we have \ $U(s,t)(\omega)\to 0$ \ as \ $(s,t)\to (1,t_0)$ \ with \ $t_0 \in [0,1)$ \ for all
 \ $\omega\in C_3$, \ where \ $C_3$ \ denotes the set of those \ $\omega\in\Omega$ \ for which \eqref{help_LIL_Wiener_field1} holds
 (then \ $\PP(C_3)=1$).
\ Similarly, if \ $(s,t)\to (s_0,t_0)$ \ with \ $s_0 \in[0,1)$ \ and \ $t_0 = 1$, \ then
 \ $U(s,t)(\omega)\to 0$ \ for all \ $\omega\in C_3$.
\ To conclude, note that
 \[
   C_1\cap C_2\cap C_3 \subseteq \{ \omega\in\Omega : [0,1]^2\ni(s,t)\mapsto U(s,t)(\omega) \;\;\text{is continuous} \},
 \]
 and \ $\PP(C_1\cap C_2\cap C_3)=1$, \ yielding that \ $U$ \ has continuous sample paths almost surely.
\proofend
\end{Rem}

\begin{Rem}\label{Rem_counter}
In what follows we show that the form of the space-domain transformations used for representing
 a tied-down Wiener bridge as a space-domain scaled stationary Ornstein-Uhlenbeck field is not appropriate
 in case of a bivariate Wiener bridge.
More precisely, one cannot find functions \ $f,g:(0,1]\to(0,\infty)$ \ such that
 \ $f$ \ is monotone and the random field
\begin{align*}
  V(s,t)
    :=\begin{cases}
        \sqrt{g(s)g(t)} X\left(\ln(f(s)), \ln(f(t)) \right),
              & \text{if \ $(s,t)\in(0,1]^2$,}\\
            0 & \text{if \ $s = 0$ \ or \ $t = 0$,}
      \end{cases}
 \end{align*}
 generates the same law on \ $C([0,1]^2)$ \ as a bivariate Wiener bridge \ $B$,
 \ where \ $X$ \ is a stationary Ornstein-Uhlenbeck field with parameters \ $(\frac{1}{2},\frac{1}{2},1)$.
\ On the contrary, let us suppose that there exist such functions.
Without loss of generality, we may suppose that \ $f$ \ is monotone increasing.
Then, due to the covariance structure of \ $X$, \ for all \ $0 < s_1\leq s_2\leq 1$ \ and \ $0 < t_1\leq t_2\leq 1$,
 \ we have
 \begin{align*}
  s_1t_1 - s_1s_2t_1t_2
    & = \cov(B(s_1,t_1),B(s_2,t_2))
     = \cov(V(s_1,t_1),V(s_2,t_2)) \\
    & = \sqrt{g(s_1)g(s_2)g(t_1)g(t_2)} \\
    &\phantom{=}\times
         \exp\Big\{-\frac{1}{2}(\ln(f(s_2)) - \ln(f(s_1)))
                 -\frac{1}{2}(\ln(f(t_2)) - \ln(f(t_1)))\Big\},
 \end{align*}
 and hence
 \begin{align*}
 s_1t_1(1-s_2t_2)
    = \sqrt{g(s_1)g(s_2)g(t_1)g(t_2)}
        \sqrt{\frac{f(s_1)}{f(s_2)}\frac{f(t_1)}{f(t_2)} }
    =:F(s_1)G(s_2)F(t_1)G(t_2),
  \end{align*}
 with \ $F(s):=\sqrt{f(s)g(s)}$, $s\in(0,1]$, \ and \ $G(s):=\sqrt{g(s)/f(s)}$, $s\in(0,1]$.
\ Then for all \ $0< s_1\leq s_2\leq 1$ \ and \ $0< t_1\leq t_2\leq 1$, \ we have
 \[
   1 - s_2t_2 = \frac{F(s_1)}{s_1} G(s_2) \frac{F(t_1)}{t_1} G(t_2)
              =:  \widetilde F(s_1)G(s_2)\widetilde F(t_1)G(t_2).
 \]
By substituting \ $s_1=t_1=t_2:=\frac{1}{2}$, \ and \ \ $s_1=t_1:=\frac{1}{2}$, \ $t_2:=1$, \ we have
 \[
    1- \frac{s_2}{2} = (\widetilde F(1/2))^2 G(1/2) G(s_2),
    \qquad s_2\in[1/2,1],
 \]
 and
  \[
    1 - s_2 = (\widetilde F(1/2))^2 G(1) G(s_2),
    \qquad s_2\in[1/2,1],
 \]
 respectively.
Consequently, \ $(\widetilde F(1/2))^2 G(1/2)\ne 0$ \ and \ $(\widetilde F(1/2))^2 G(1)\ne 0$,
 \[
   G(s_2) = \frac{1-s_2/2}{(\widetilde F(1/2))^2 G(1/2)}
   \qquad \text{and}\qquad
   G(s_2) = \frac{1-s_2}{(\widetilde F(1/2))^2 G(1)},
   \qquad s_2\in[1/2,1],
 \]
 which yields us to a contradiction (by choosing, e.g., \ $s_2=1$).
\proofend
\end{Rem}

\section{Tied-down scaled Wiener bridges}

Let \ $S>0$, $T>0$, \ and \ $\alpha>0$, $\beta>0$, \ and let us consider
 a zero-mean Gauss field \ $\{X^{(\alpha,\beta)}(s,t) : (s,t)\in[0,S]\times[0,T]\}$ \
 with continuous sample paths almost surely and with covariance function
 \[
  \cov(X^{(\alpha,\beta)}(s_1,t_1), X^{(\alpha,\beta)}(s_2,t_2))
     = R^{(\alpha)}_S(s_1,s_2)R^{(\beta)}_T(t_1,t_2)
 \]
 for \ $(s_1,t_1), (s_2,t_2)\in[0,S]\times[0,T]$,
 \ where \ $R^{(\alpha)}_S$ \ is the covariance function of a scaled Wiener bridge \ $X^{(\alpha)}$ \ on \ $[0,S]$ \ with parameter
  \ $\alpha$ \ given by
 \[
   R^{(\alpha)}_S(s_1,s_2)
   = \begin{cases}
      \frac{(S-s_1)^\alpha (S-s_2)^\alpha}{1-2\alpha}
           (S^{1-2\alpha} - (S-(s_1\wedge s_2))^{1-2\alpha})
           & \text{if \ $\alpha\ne \frac{1}{2}$,}\\[1mm]
        \sqrt{(S-s_1)(S-s_2)}\ln\left(\frac{S}{S-(s_1\wedge s_2)}\right)
           & \text{if \ $\alpha =\frac{1}{2}$,}
     \end{cases}
 \]
 for \ $s_1,s_2\in[0,S]$, \ and \ $R^{(\beta)}_T$ \ is the covariance function of a scaled Wiener bridge \ $X^{(\beta)}$ \ on \ $[0,T]$ \
 with parameter \ $\beta$.
\ Here \ $R^{(\alpha)}_S$ \ is defined to be \ $0$ \ on the line segments between \ $(0,S)$ \ and \ $(S,S)$,
 \ and \ $(S,0)$ \ and \ $(S,S)$, \ respectively, as a consequence of \ $\lim_{(s_1,s_2)\to (s,S)} R^{(\alpha)}_S(s_1,s_2)
  = \lim_{(s_1,s_2)\to (S,s)} R^{(\alpha)}_S(s_1,s_2) = 0$, \ $s\in[0,S]$ \ (for a detailed discussion, see Barczy and Igl\'oi \cite{BarIgl}).
We note that scaled Wiener bridges were introduced by Brennan and Schwartz \cite{BreSch}, and see also Mansuy \cite{Man};
 and the random field \ $X^{(\alpha,\beta)}$ \ has already been introduced in Barczy and Igl\'oi \cite[page 5]{BarIgl}.
Since for independent scaled Wiener bridges \ $(X^{(\alpha)}(s))_{s\in[0,S]}$ \ and \ $(X^{(\beta)}(t))_{t\in[0,T]}$,
 \ the random (but not Gauss) field \ $\{ X^{(\alpha)}(s)X^{(\beta)}(t) : (s,t)\in[0,S]\times[0,T]\}$ \ admits the same
 covariances as \ $X^{(\alpha,\beta)}$, \ there exists a zero mean Gauss field with the given covariances.
Later on (see Proposition \ref{Prop_scaled_Wiener_field}), we will see that the continuity assumption
 can also be fulfilled.
Note that \ $X^{(\alpha,\beta)}$ \ is zero on the border of a rectangle with vertices \ $(0,0)$, $(0,S)$, $(0,T)$, \ and \ $(S,T)$,
 so we can call it a tied-down scaled Wiener bridge with parameters \ $(\alpha,\beta)$.
\ This class of Gauss processes may deserve more attention since it would generalize some well-known
 limit processes in mathematical statistics such as a Kiefer process, see, e.g.,
 Deheuvels et al. \cite[formula (3.8) with $\gamma=\delta=0$]{DehPecYor} or Shorack and Wellner \cite[Exercise 12, page 32]{ShoWel}.
In Remark \ref{Rem_Kiefer}, we detail the case of a Kiefer process.

The following result can be considered as a generalization of the corresponding one for scaled Wiener bridges
 in Subsection 3.1 in Barczy and Kern \cite{BarKer3}.

\begin{Pro}\label{Prop_scaled_Wiener_field}
Let \ $\{X(s,t) : s,t\in\RR\}$ \ be a stationary Ornstein-Uhlenbeck field with parameters \ $(\frac{1}{2},\frac{1}{2},1)$ \ represented as in Proposition \ref{prop_stac_OU_field}.
Let \ $S>0$, $T>0$, \ and \ $\alpha >0$, $\beta>0$.
Then the random field
 \begin{align*}
  U(s,t)
    :=\begin{cases}
        \sqrt{g^{(\alpha)}_S(s) g^{(\beta)}_T(t)}\, X\left(\ln\big(f^{(\alpha)}_S(s)\big) , \ln\big(f^{(\beta)}_T(t)\big) \right),
              & \text{if \ $(s,t)\in(0,S)\times (0,T)$,}\\
            0 & \text{if \ $s\in\{0,S\}$ \ or \ $t\in\{0,T\}$,}
      \end{cases}
 \end{align*}
 is a tied-down scaled Wiener bridge with parameters \ $(\alpha,\beta)$, \ where
  \[
  \sqrt{g^{(\alpha)}_S(s)}
     := \begin{cases}
            (S-s)^\alpha \sqrt{\frac{S^{1-2\alpha} - (S-s)^{1-2\alpha}}{1-2\alpha}}
               & \text{if \ $\alpha\ne\frac{1}{2}$,}\\[1mm]
           \sqrt{(S-s)\ln\left(\frac{S}{S-s}\right)}
               & \text{if \ $\alpha=\frac{1}{2}$,}
        \end{cases}
       \qquad s\in(0,S),
  \]
  \[
    f^{(\alpha)}_S(s)
       :=\begin{cases}
           \frac{S^{2\alpha}}{1-2\alpha}(S^{1-2\alpha} - (S-s)^{1-2\alpha})
                & \text{if \ $\alpha\ne\frac{1}{2}$,}\\
        S\ln\left(\frac{S}{S-s}\right)
               & \text{if \ $\alpha=\frac{1}{2}$,}
        \end{cases}
       \qquad s\in(0,S),
  \]
  and \ $\sqrt{g^{(\beta)}_T(t)}$ \ and \ $f^{(\beta)}_T(t)$ \ are defined similarly as \ $\sqrt{g^{(\alpha)}_S(s)}$ \
  and \ $f^{(\alpha)}_S(s)$, \ respectively, by replacing \ $\alpha$ \ by \ $\beta$, \
   $s$ \ by \ $t$, \ and \ $S$ \ by \ $T$.
\end{Pro}

\noindent{\bf Proof.}
First, we check that both fields \ $U$ \ and \ $X^{(\alpha,\beta)}$ \ are zero-mean Gauss fields on \ $[0,S]\times[0,T]$ \ with
 continuous sample paths almost surely.
The only property that does not follow immediately is that \ $U$ \ has continuous sample paths almost surely.
For all \ $(s,t)\in(0,S)\times(0,T)$,
 \begin{align}\label{help1}
  \begin{split}
    U(s,t)
     &=\sqrt{g^{(\alpha)}_S(s) g^{(\beta)}_T(t)}\,
        X\left(\ln\left( f^{(\alpha)}_S(s)\right), \ln\left(f^{(\beta)}_T(t)\right) \right) \\
    & = \sqrt{g^{(\alpha)}_S(s) g^{(\beta)}_T(t)}
        \exp\Big\{ -\frac{1}{2}\ln\left( f^{(\alpha)}_S(s) \right) - \frac{1}{2}\ln\left( f^{(\beta)}_T(t) \right) \Big\}
        W\left(\ee^{\ln\left( f^{(\alpha)}_S(s) \right)},\ee^{\ln\left(f^{(\beta)}_T(t)\right)}\right) \\
     & = \sqrt{g^{(\alpha)}_S(s) g^{(\beta)}_T(t)}
         \frac{1}{\sqrt{f^{(\alpha)}_S(s) f^{(\beta)}_T(t)}}
          W\left( f^{(\alpha)}_S(s) , f^{(\beta)}_T(t) \right) \\
     & = \left(1-\frac{s}{S}\right)^\alpha
         \left(1-\frac{t}{T}\right)^\beta
          W\left( f^{(\alpha)}_S(s) , f^{(\beta)}_T(t) \right),
  \end{split}
 \end{align}
 and, by an easy calculation,
 \begin{align}\label{help3}
  \lim_{s\downarrow 0} f^{(\alpha)}_S(s)
     = 0, \quad \alpha>0,
    \qquad \text{and}\qquad
   \lim_{s\uparrow S} f^{(\alpha)}_S(s)
     = \begin{cases}
         \frac{S}{1-2\alpha} & \text{if \ $0<\alpha<\frac{1}{2}$,}\\
         +\infty & \text{if \ $\alpha\geq \frac{1}{2}$,}
       \end{cases}
 \end{align}
 and similar expressions hold for \ $\lim_{t\downarrow 0} f^{(\beta)}_T(t)$ \ and \ $\lim_{t\uparrow T} f^{(\beta)}_T(t)$.
In what follows, we will proceed similarly as in Remark \ref{Rem_alternativ_proof}.
From Section \ref{Section2} recall the notations
 \begin{align*}
   &C_1=\{ \omega\in\Omega : \RR_+^2\ni(s,t)\mapsto W(s,t)(\omega)\;\;\text{is continuous}\},\\
   &C_2=\{ \omega\in\Omega : \text{\eqref{help_LIL_Wiener_field2} holds} \},
       \qquad C_3=\{ \omega\in\Omega : \text{\eqref{help_LIL_Wiener_field1} holds} \},
 \end{align*}
 and we have \ $\PP(C_1) = \PP(C_2) = \PP(C_3) = 1$.

\noindent If \ $(s,t)\to (s_0,t_0)\in[0,S]\times [0,T]$ \ with \ $s_0\ne S$ \ and \ $t_0\ne T$
 \ and \ $\alpha>0$ \ and \ $\beta>0$,
 \ then, by \eqref{help1}, we have \ $U(s,t)(\omega)\to U(s_0,t_0)(\omega)$ \ for all \ $\omega\in C_1$.

\noindent If \ $(s,t)\to (S,T)$ \ and \ $0<\alpha<\frac{1}{2}$ \ and \ $0<\beta<\frac{1}{2}$,
 \ then, by \eqref{help1} and \eqref{help3}, we have \ $U(s,t)(\omega) \to 0\cdot W(S/(1-2\alpha), T/(1-2\beta))(\omega) = 0$ \
 for all \ $\omega\in C_1$.

\noindent If \ $(s,t)\to (S,T)$ \ and \ $\alpha\geq \frac{1}{2}$ \ and \ $\beta\geq \frac{1}{2}$,
 \ then, by \eqref{help1} and \eqref{help_LIL_Wiener_field2}, we have
 \begin{align*}
 U(s,t)(\omega)
   &= \frac{2}{S^\alpha T^\beta}
      \sqrt{(S-s)^{2\alpha} (T-t)^{2\beta}  f^{(\alpha)}_S(s) f^{(\beta)}_T(t) \ln(\ln(f^{(\alpha)}_S(s) f^{(\beta)}_T(t)))} \\
   &\phantom{=\;}
      \times\frac{W(f^{(\alpha)}_S(s), f^{(\beta)}_T(t))(\omega)}
           {\sqrt{4 f^{(\alpha)}_S(s) f^{(\beta)}_T(t) \ln(\ln(f^{(\alpha)}_S(s) f^{(\beta)}_T(t)))} }
     \to 0
 \end{align*}
 for all \ $\omega\in C_2$, \ since
 \[
   \lim_{s\uparrow S,\, t\uparrow T}
    (S-s)^{2\alpha} (T-t)^{2\beta}  f^{(\alpha)}_S(s) f^{(\beta)}_T(t) \ln(\ln(f^{(\alpha)}_S(s) f^{(\beta)}_T(t)))
    =0,\qquad \alpha\geq \frac{1}{2}, \quad \beta\geq \frac{1}{2}.
  \]
Indeed, in case of \ $\alpha>\frac{1}{2}$ \ and \ $\beta>\frac{1}{2}$,
 \begin{align*}
  &\lim_{s\uparrow S, \,t\uparrow T}
     (S-s)^{2\alpha} (T-t)^{2\beta}  f^{(\alpha)}_S(s) f^{(\beta)}_T(t) \ln(\ln(f^{(\alpha)}_S(s) f^{(\beta)}_T(t))) \\
  & = \lim_{s\uparrow S,\, t\uparrow T}
     \frac{S^{2\alpha} T^{2\beta}}{(2\alpha - 1)(2\beta - 1)}
     (S-s)^{2\alpha} (T-t)^{2\beta}
     ((S-s)^{1-2\alpha} - S^{1-2\alpha})^{\frac{2\alpha}{2\alpha-1}} ((T-t)^{1-2\beta} - T^{1-2\beta})^{\frac{2\beta}{2\beta-1}}\\
  &\phantom{= \lim_{s\uparrow S,\, t\uparrow T}\;}
     \times\frac{\ln \left(\ln \left( \frac{S^{2\alpha} T^{2\beta}}{(2\alpha - 1)(2\beta - 1)}
                                ((S-s)^{1-2\alpha} - S^{1-2\alpha})
                                ((T-t)^{1-2\beta} - T^{1-2\beta})
               \right) \right)}
          {((S-s)^{1-2\alpha} - S^{1-2\alpha})^{\frac{1}{2\alpha-1}} ((T-t)^{1-2\beta} - T^{1-2\beta})^{\frac{1}{2\beta-1}}}
   = 0,
 \end{align*}
 where we used that
 \begin{align*}
  & \lim_{s\uparrow S}
     (S-s)^{2\alpha}((S-s)^{1-2\alpha} - S^{1-2\alpha})^{\frac{2\alpha}{2\alpha-1}}
   = \lim_{s\uparrow S}
     \left(1- S^{1-2\alpha}(S-s)^{2\alpha - 1} \right)^{\frac{2\alpha}{2\alpha-1}}
   = 1, \\
 &\lim_{t\uparrow T}
     (T-t)^{2\beta}((T-t)^{1-2\beta} - T^{1-2\beta})^{\frac{2\beta}{2\beta-1}}
   = 1,
 \end{align*}
 and
 \[
  \lim_{h\to\infty} \frac{1}{h^\varepsilon}\ln(\ln(h)) = 0, \qquad \forall \; \varepsilon>0;
 \]
 in case of \ $\alpha=\frac{1}{2}$ \ and \ $\beta=\frac{1}{2}$, \ by L'Hospital's rule,
 \begin{align*}
  &\lim_{s\uparrow S,\, t\uparrow T}
     (S-s)(T-t)f^{(1/2)}_S(s) f^{(1/2)}_T(t) \ln(\ln(f^{(1/2)}_S(s) f^{(1/2)}_T(t))) \\
  & = \lim_{s\uparrow S,\, t\uparrow T}
     ST (S-s)(T-t)\left( \ln\left(\frac{S}{S-s}\right)
                         \ln\left(\frac{T}{T-t}\right)
                  \right)^2 \\
  &\phantom{= \lim_{s\uparrow S, \,t\uparrow T}\;}
     \times\frac{\ln \left(\ln \left( ST\ln\left(\frac{S}{S-s}\right)
                         \ln\left(\frac{T}{T-t}\right)
               \right) \right)}
          { \ln\left(\frac{S}{S-s}\right) \ln\left(\frac{T}{T-t}\right) }
   = 0,
 \end{align*}
 and the other cases can be handled similarly.

\noindent If \ $(s,t)\to (S,t_0)$ \ with \ $t_0\in[0,T)$ \ and \ $0<\alpha<\frac{1}{2}$, \ $\beta>0$, \ then,
 by \eqref{help1} and \eqref{help3}, we have \ $U(s,t)(\omega) \to 0(1-t_0/\beta)^\beta W(S/(1-2\alpha), f^{(\beta)}_T(t_0))(\omega) = 0$ \
 for all \ $\omega\in C_1$.

\noindent If \ $(s,t)\to (S,t_0)$ \ with \ $t_0\in[0,T)$ \ and \ $\alpha\geq \frac{1}{2}$, \ $\beta>0$, \ then,
 by \eqref{help1}, we have
 \begin{align*}
  U(s,t)
    &= \left(1 - \frac{s}{S}\right)^\alpha
      \left(1 - \frac{t}{T}\right)^\beta
       \sqrt{4f^{(\alpha)}_S(s)(f^{(\beta)}_T(t_0) + 1)\ln\big(\ln(f^{(\alpha)}_S(s) (f^{(\beta)}_T(t_0) + 1))\big)} \\
    &\phantom{=\;}
      \times \frac{W(f^{(\alpha)}_S(s),f^{(\beta)}_T(t))}
        {\sqrt{4f^{(\alpha)}_S(s)(f^{(\beta)}_T(t_0) + 1)\ln\big(\ln(f^{(\alpha)}_S(s) (f^{(\beta)}_T(t_0) + 1) )\big)} },
 \end{align*}
 where
 \[
    (f^{(\alpha)}_S(s),f^{(\beta)}_T(t))
       \in D^{(1)}_{f^{(\alpha)}_S(s)(f^{(\beta)}_T(t_0) + 1)}
 \]
 provided that \ $s$ \ is sufficiently close to \ $S$ \ (it is enough to choose \ $s$ \ such that \ $f^{(\alpha)}_S(s)\geq 1$ \ which
 can be done due to \eqref{help3}) and \ $t$ \ is sufficiently close to \ $t_0$.
By \eqref{help_LIL_Wiener_field1}, using the calculations for the case \ $(s,t)\to (S,T)$ \ and $\alpha\geq \frac{1}{2}$,
 $\beta\geq\frac{1}{2}$, \ as well, we have \ $U(s,t)(\omega)\to 0$ \ as \ $(s,t)\to (S,t_0)$ \ with \ $t_0\in[0,T)$ \ for all \ $\omega\in C_3$.

\noindent Similarly, if \ $(s,t)\to (s_0,T)$ \ with \ $s_0\in[0,S)$, \ then \ $U(s,t)(\omega)\to 0$ \ for all \ $\omega\in C_3$.
Since
 \[
    C_1\cap C_2 \cap C_3
       \subseteq \{ \omega\in\Omega : [0,S]\times [0,T]\ni (s,t)\mapsto U(s,t)(\omega)\;\;\text{is continuous}\} ,
 \]
 and \ $\PP(C_1\cap C_2\cap C_3)=1$, \ we have \ $U$ \ has continuous sample paths almost surely.

To conclude, it is sufficient to check that the covariance functions of \ $U$ \ and \ $X^{(\alpha,\beta)}$ \ coincide.
First let us suppose that \ $\alpha\ne\frac{1}{2}$ \ and \ $\beta\ne\frac{1}{2}$.
\ Then for all \ $0< s_1\leq s_2<S$ \ and \ $0<t_1\leq t_2< T$ \ (which can be supposed without loss of generality),
 we have
 \begin{align*}
   &\cov(U(s_1,t_1),U(s_2,t_2))\\
   & =(S-s_1)^\alpha \sqrt{\frac{S^{1-2\alpha} - (S-s_1)^{1-2\alpha}}{1-2\alpha}}
         (S-s_2)^\alpha \sqrt{\frac{S^{1-2\alpha} - (S-s_2)^{1-2\alpha}}{1-2\alpha}} \\
   &\phantom{=\,}
       \times (T-t_1)^\beta \sqrt{\frac{T^{1-2\beta} - (T-t_1)^{1-2\beta}}{1-2\beta}}
        (T-t_2)^\beta \sqrt{\frac{T^{1-2\beta} - (T-t_2)^{1-2\beta}}{1-2\beta}}\\
 &\phantom{=\,}
       \times \exp\left\{ -\frac{1}{2}
            \ln\left( \frac{S^{1-2\alpha} - (S-s_2)^{1-2\alpha}}{S^{1-2\alpha} - (S-s_1)^{1-2\alpha}}\right)
            -\frac{1}{2}
            \ln\left( \frac{T^{1-2\beta} - (T-t_2)^{1-2\beta}}{T^{1-2\beta} - (T-t_1)^{1-2\beta}}\right)
           \right\}
 \end{align*}
 \begin{align*}
 & = \frac{(S-s_1)^{\alpha} (S-s_2)^{\alpha}}{1-2\alpha}
      \frac{(T-t_1)^{\beta} (T-t_2)^{\beta}}{1-2\beta}
      (S^{1-2\alpha} - (S-s_1)^{1-2\alpha})
      (T^{1-2\beta} - (T-t_1)^{1-2\beta})\\
  & = \cov(X^{(\alpha,\beta)}(s_1,t_1), X^{(\alpha,\beta)}(s_2,t_2)),
 \end{align*}
 as desired, where we used that \ $f^{(\alpha)}_S$ \ and \ $f^{(\beta)}_T$ \ are strictly increasing.
Let us suppose now that \ $\alpha=\beta=\frac{1}{2}$.
\ Then for all \ $0< s_1\leq s_2< S$ \ and \ $0< t_1\leq t_2< T$ \ (which can be supposed without loss of generality),
 we have
 \begin{align*}
   \cov(U(s_1,t_1),U(s_2,t_2))
   &= \sqrt{(S-s_1)\ln\left(\frac{S}{S-s_1}\right)}
     \sqrt{(S-s_2)\ln\left(\frac{S}{S-s_2}\right)}\\
   &\phantom{=\,}
     \times\sqrt{(T-t_1)\ln\left(\frac{T}{T-t_1}\right)}
     \sqrt{(T-t_2)\ln\left(\frac{T}{T-t_2}\right)}\\
    &\phantom{=\,}
     \times\exp\left\{ -\frac{1}{2}\ln\left(\frac{\ln\left(\frac{S}{S-s_2}\right)}{\ln\left(\frac{S}{S-s_1}\right)}\right)
                 -\frac{1}{2}\ln\left(\frac{\ln\left(\frac{T}{T-t_2}\right)}{\ln\left(\frac{T}{T-t_1}\right)}\right)
                \right\} \\
  &= \sqrt{(S-s_1)(S-s_2)(T-t_1)(T-t_2)}
      \ln\left(\frac{S}{S-s_1}\right)
      \ln\left(\frac{T}{T-t_1}\right)\\
  & = \cov(X^{(\alpha,\beta)}(s_1,t_1), X^{(\alpha,\beta)}(s_2,t_2)),
 \end{align*}
 as desired.
The cases \ $\alpha\ne\frac{1}{2}$, \ $\beta=\frac{1}{2}$, \ and \ $\alpha=\frac{1}{2}$, \ $\beta\ne\frac{1}{2}$, \
 can be handled similarly.
\proofend

\begin{Rem}
Note that if \ $\alpha=1$, $\beta=1$, \ and \ $S=T=1$, \ then
 \ $g^{(\alpha)}_S(s) = s(1-s)$, $s\in(0,1)$, \ $g^{(\beta)}_T(t) = t(1-t)$, \ $t\in(0,1)$,
 \ and \ $f^{(\alpha)}_S(s) = \frac{s}{1-s}$, $s\in(0,1)$, \
 \ $f^{(\beta)}_T(t) = \frac{t}{1-t}$, $t\in(0,1)$.
Hence, in case of  \ $\alpha=1$, $\beta=1$, \ and \ $S=T=1$, \
 Proposition \ref{Prop_scaled_Wiener_field} gives back Proposition \ref{prop_tied_dow_W-bridge} (as expected).
\proofend
\end{Rem}

The next remark is devoted to the case of a Kiefer process.

\begin{Rem}\label{Rem_Kiefer}
Let \ $\{W(s,t) : s,t\in\RR_+\}$ \ be a standard Wiener field.
Then the random field
 \[
   \Big\{ K(s,t):= W(s,t) - sW(1,t) : s\in[0,1],\, t\in\RR_+ \Big\},
 \]
 is a zero-mean Gauss field with continuous sample paths almost surely and with covariance function
 \[
   \cov(K(s_1,t_1), K(s_2,t_2))
     =(s_1\wedge s_2 - s_1s_2)(t_1\wedge t_2),
    \qquad s_1,s_2\in[0,1],\, t_1,t_2\in\RR_+.
 \]
The random field \ $K$ \ is known as a Kiefer process, see, e.g., Deheuvels et al. \cite[formula (3.8) with $\gamma=\delta=0$]{DehPecYor}
 or Shorack and Wellner \cite[Exercise 12, page 32]{ShoWel}.
Note that, {\sl formally}, with \ $S=1$, $T=\infty$, \ $\alpha=1$ \ and \ $\beta=0$, \ we have
 \ $\cov(K(s_1,t_1), K(s_2,t_2)) = R^{(\alpha)}_1(s_1,s_2) R^{(\beta)}_\infty(t_1,t_2)$, \ $s_1,s_2\in[0,1]$, $t_1,t_2\in\RR_+$.
\ Further, similarly as in the proof of Proposition \ref{Prop_scaled_Wiener_field},
 one can check that the random field
 \begin{align*}
  U(s,t)
    :=\begin{cases}
        \sqrt{s(1-s)t} \,X\left(\ln\left(\frac{s}{1-s}\right) , \ln(t) \right),
              & \text{if \ $(s,t)\in(0,1)\times (0,\infty)$,}\\
            0 & \text{if \ $s\in\{0,1\}$ \ or \ $t=0$,}
      \end{cases}
 \end{align*}
 generates the same law on \ $C([0,1]\times[0,\infty))$ \ as a Kiefer process.
Indeed, both fields \ $U$ \ and \ $K$ \ are zero-mean Gauss fields on \ $[0,1]\times[0,\infty)$ \ with
 continuous sample paths almost surely (which can be checked similarly as in the proof of Proposition
 \ref{Prop_scaled_Wiener_field}), and for all \ $0< s_1\leq s_2< 1$ \ and \ $0< t_1\leq t_2< \infty$ \
 (which can be supposed without loss of generality), we have
 \begin{align*}
  \cov(U(s_1,t_1),U(s_2,t_2))
    & = \sqrt{s_1(1-s_1)t_1}\sqrt{s_2(1-s_2)t_2}\\
    &\phantom{=\,}
        \times\exp\left\{ -\frac{1}{2}\left(\ln\left(\frac{s_2}{1-s_2}\right)- \ln\left(\frac{s_1}{1-s_1}\right)\right)
                 -\frac{1}{2}(\ln(t_2) - \ln(t_1)) \right\} \\
    & = \sqrt{s_1t_1s_2t_2(1-s_1)(1-s_2)}
       \sqrt{\frac{s_1}{1-s_1} \frac{1-s_2}{s_2} \frac{t_1}{t_2}}\\
    & = s_1t_1(1-s_2)
      = \cov(K(s_1,t_1),K(s_2,t_2)),
 \end{align*}
 as desired.
Note that, {\sl formally}, this result is nothing else but Proposition \ref{Prop_scaled_Wiener_field}
 with \ $S=1$, $T=\infty$, \ $\alpha=1$ \ and \ $\beta=0$.
\proofend
\end{Rem}

\section{$(F,G)$-Wiener bridges}
Let \ $f:\RR_+\to\RR_+$ \ and \ $g:\RR_+\to\RR_+$ \ be probability density functions on \ $\RR_+$ \ and
 let us consider the corresponding cumulative distribution functions \ $F:\RR_+\to[0,1]$, $F(s):=\int_0^s f(u)\,\dd u$,
 $s\in\RR_+$, \ and \ $G:\RR_+\to[0,1]$, $G(t):=\int_0^t g(u)\,\dd u$, $t\in\RR_+$.
Further, let
 \[
   S:=\inf\{s\in\RR_+ : F(s)=1\}\in(0,\infty],
     \quad \quad
   T:=\inf\{t\in\RR_+ : G(t)=1\}\in(0,\infty],
 \]
 with the convention \ $\inf\emptyset:=\infty$.
\ Let us assume that \ $f$ \ and \ $g$ \ are continuous on \ $[0,T)$ \ and \ $[0,S)$, \ respectively,
 and that there exist \ $\delta_1\in(0,S)$ \ and \ $\delta_2\in(0,T)$
 \ such that \ $f(t)\ne 0$ \ for all \ $t\in(0,\delta_1)$, \ and \ $g(t)\ne 0$ \ for all \ $t\in(0,\delta_2)$.

Let us consider a zero-mean Gauss field \ $\{ X^{(F,G)}(s,t) : (s,t)\in[0,S)\times [0,T)\}$ \ with continuous sample
 paths almost surely and with covariance function
 \[
    \cov(X^{(F,G)}(s_1,t_1),X^{(F,G)}(s_2,t_2))
       := (F(s_1\wedge s_2) - F(s_1)F(s_2))(G(t_1\wedge t_2) - G(t_1)G(t_2))
 \]
 for \ $(s_i,t_i)\in[0,S)\times [0,T)$, $i=1,2$, \ which we call an \ $(F,G)$-Wiener bridge.
Next we check that for independent \ $F$- and \ $G$-Wiener bridges \ $(Y^{(F)}_s)_{s\in[0,S)}$ \ and
 \ $(Z^{(G)}_t)_{t\in[0,T)}$, \ the (non-Gauss) random field \ $\{Y^{(F)}_sZ^{(G)}_t : (s,t)\in[0,S)\times [0,T)\}$ \ admits the same
 covariances as \ $X^{(F,G)}$, \ and hence there exists a zero-mean Gauss field with the given covariances.
For the existence and properties of an \ $F$-Wiener bridge \ $(Y^{(F)}_s)_{s\in[0,S)}$ \ under the given conditions on \ $f$, \ see Subsection 3.3
 in Barczy and Kern \cite{BarKer3}, Shorack and Wellner \cite[page 838]{ShoWel}, van der Vaart \cite[page 266]{Vaa} or Khmaladze \cite[equation (4)]{Khm}.
Here we only recall that \ $(Y^{(F)}_s)_{s\in[0,S)}$ \ is a zero-mean Gauss process having continuous sample paths almost surely
 and covariance function \ $F(s\wedge t) - F(s)F(t)$, $s,t\in[0,S)$, \ satisfying \ $Y^{(F)}_0=0$ \ and
 \ $\PP(\lim_{s\uparrow S} Y^{(F)}_s = 0)=1$.
\ Then for \ $(s_i,t_i)\in[0,S)\times [0,T)$, $i=1,2$, \ we have
 \begin{align*}
  \cov(Y^{(F)}_{s_1}Z^{(G)}_{t_1},Y^{(F)}_{s_2}Z^{(G)}_{t_2})
    & = \EE(Y^{(F)}_{s_1}Z^{(G)}_{t_1}Y^{(F)}_{s_2}Z^{(G)}_{t_2})
      = \EE(Y^{(F)}_{s_1}Y^{(F)}_{s_2})\EE(Z^{(G)}_{t_1}Z^{(G)}_{t_2})\\
    & = (F(s_1\wedge s_2) - F(s_1)F(s_2))
       (G(t_1\wedge t_2) - G(t_1)G(t_2)) .
 \end{align*}
We will see that the continuity assumption on the sample paths of \ $X^{(F,G)}$ \ can also be fulfilled,
 and we will give a possible motivation of the name \ $(F,G)$-Wiener bridge as well,
 see Proposition \ref{prop_FG_Wiener_field} and paragraph just after it, respectively.

The following result can be considered as a generalization of  the corresponding one for $F$-Wiener bridges
 in Subsection 3.3 in Barczy and Kern \cite{BarKer3}.

\begin{Pro}\label{prop_FG_Wiener_field}
Let \ $\{X(s,t) : s,t\in\RR\}$ \ be a stationary Ornstein-Uhlenbeck field with parameters
 \ $(\frac{1}{2},\frac{1}{2},1)$ \ represented as in Proposition \ref{prop_stac_OU_field}.
Then the random field
 \begin{align*}
  U(s,t)
    := \sqrt{F(s)(1-F(s))G(t)(1-G(t))}
        X\left(\ln\left(\frac{F(s)}{1-F(s)}\right) , \ln\left(\frac{G(t)}{1-G(t)}\right) \right)
 \end{align*}
 for \ $(s,t)\in(0,S)\times (0,T)$, \ and \ $U(s,t):=0$ \ for \ $s=0$ \ or \ $t=0$ \ is an \ $(F,G)$-Wiener bridge.
\end{Pro}

\noindent{\bf Proof.}
First, we check that both fields \ $U$ \ and \ $X^{(F,G)}$ \ are zero-mean Gauss fields on \ $[0,S)\times [0,T)$ \
 with continuous sample paths almost surely.
The only property that does not follow immediately is that \ $U$ \ has continuous sample paths almost surely.
For all \ $(s,t)\in(0,S)\times (0,T)$,
 \begin{align}\label{help4}
  \begin{split}
    U(s,t)
     & =\sqrt{F(s)(1-F(s))G(t)(1-G(t))}
        X\left(\ln\left(\frac{F(s)}{1-F(s)}\right) , \ln\left(\frac{G(t)}{1-G(t)}\right) \right) \\
    & = \sqrt{F(s)(1-F(s))G(t)(1-G(t))}
        \exp\Big\{ -\frac{1}{2}\ln\left(\frac{F(s)}{1-F(s)}\right) - \frac{1}{2}\ln\left(\frac{G(t)}{1-G(t)}\right) \Big\}\\
    &\phantom{=\;}
        \times W\left(\ee^{\ln\left(\frac{F(s)}{1-F(s)}\right)},\ee^{\ln\left(\frac{G(t)}{1-G(t)}\right)}\right) \\
     & = \sqrt{F(s)(1-F(s))G(t)(1-G(t))}
        \sqrt{\frac{1-F(s)}{F(s)}\cdot \frac{1-G(t)}{G(t)} }\,
        W\left(\frac{F(s)}{1-F(s)},\frac{G(t)}{1-G(t)}\right) \\
     & = (1-F(s))(1-G(t)) W\left(\frac{F(s)}{1-F(s)},\frac{G(t)}{1-G(t)}\right) .
  \end{split}
 \end{align}
Since \ $F$ \ and \ $G$ \ are continuous, \ $F(0)=G(0)=0$, \ we get
 \begin{align*}
  C_1 & = \{ \omega\in\Omega : \RR_+^2\ni(s,t)\mapsto W(s,t)(\omega)\;\;\text{is continuous}\} \\
      & \subseteq \{ \omega\in\Omega : [0,S)\times[0,T)\ni(s,t)\mapsto U(s,t)(\omega)\;\;\text{is continuous}\},
 \end{align*}
 and consequently, due to \ $\PP(C_1)=1$, \ the sample paths of \ $U$ \ are continuous almost surely.
To conclude, it is enough to check that the covariance functions of \ $U$ \ and \ $X^{(F,G)}$ \ coincide.
For all \ $0 < s_1\leq s_2<S$ \ and \ $0 < t_1\leq t_2<T$ \ (which can be assumed without loss of generality),
 we have
 \begin{align*}
 \cov(U(s_1,t_1),U(s_2,t_2))
    &= \sqrt{F(s_1)(1-F(s_1))G(t_1)(1-G(t_1))}\\
    &\phantom{=}\times \sqrt{F(s_2)(1-F(s_2))G(t_2)(1-G(t_2))}\\
    &\phantom{=}\times
    \exp\left\{ -\frac{1}{2}\ln\left(\frac{F(s_2)(1-F(s_1))}{F(s_1)(1-F(s_2))}\right)
                -\frac{1}{2}\ln\left(\frac{G(t_2)(1-G(t_1))}{G(t_1)(1-G(t_2))}\right) \right\}\\
    & = F(s_1)(1-F(s_2))G(t_1)(1-G(t_2)),
 \end{align*}
 as desired, where we used that the functions \ $(0,S)\ni s \mapsto \frac{F(s)}{1-F(s)}$ \ and
 \ $(0,T)\ni t \mapsto \frac{G(t)}{1-G(t)}$ \ are monotone increasing.
\proofend

Concerning the name \ $(F,G)$-Wiener bridge for \ $\{ X^{(F,G)}(s,t) : (s,t)\in[0,S)\times [0,T)\}$, \ we point out that
 \ $X^{(F,G)}(s,t) \to 0$ \ almost surely as \ $\frac{s}{S}\vee \frac{t}{T}\to 1$, \ which can be seen using \eqref{help4}
 and similar arguments as in Remark \ref{Rem_alternativ_proof}.

\begin{Rem}
Let \ $\{U(s,t) : s,t\in\RR_+\}$ \ be a zero-mean Gauss field with continuous sample paths almost surely
 and with covariance function of product type
 \[
   \cov(U(s_1,t_1), U(s_2,t_2)) = c(s_1,s_2)\widetilde c(t_1,t_2),\qquad s_1,s_2,t_1,t_2\in\RR_+,
 \]
 where \ $c:\RR_+\to\RR$ \ and \ $\widetilde c:\RR_+\to\RR$ \ are some (appropriately) given functions.
Similarly as in Remark \ref{Rem_counter}, one can ask whether there exist functions \ $f,\widetilde f:\RR_+\to(0,\infty)$
 \ and \ $g,\widetilde g:\RR_+\to\RR_+$ \ such that \ $f$ \ and \ $\widetilde f$ \ are monotone and the random field
\begin{align*}
  \sqrt{g(s)\widetilde g(t)} X\left(\ln(f(s)), \ln(\widetilde f(t)) \right), \qquad s,t\in\RR_+,
 \end{align*}
 generates the same law on \ $C(\RR_+^2)$ \ as \ $U$,
 \ where \ $X$ \ is a stationary Ornstein-Uhlenbeck field with parameters \ $(\frac{1}{2},\frac{1}{2},1)$.
\ Supposing that \ $f$ \ and \ $\widetilde f$ \ are monotone increasing, a necessary condition for this (following from the equality of the
 covariance functions of the random fields in question) is
 \[
   c(s_1,s_2)\widetilde c(t_1,t_2) = \sqrt{g(s_1)g(s_2)\widetilde  g(t_1) \widetilde g(t_2)}
        \sqrt{\frac{f(s_1)}{f(s_2)}\cdot \frac{\widetilde f(t_1)}{\widetilde f(t_2)} }
    =:F(s_1)G(s_2)\widetilde F(t_1)\widetilde G(t_2)
 \]
for \ $0\leq s_1\leq s_2$ \ and \ $0\leq t_1\leq t_2$, \ which can be checked similarly as in Remark \ref{Rem_counter}.
Note that all the examples presented in Propositions \ref{prop_tied_dow_W-bridge}, \ref{Prop_scaled_Wiener_field}
 and \ref{prop_FG_Wiener_field} are of this type.
\proofend
\end{Rem}

\section*{Acknowledgements}
\noindent I would like to thank Endre Igl\'oi and Peter Kern for giving useful comments on the paper.
I am undoubtedly grateful to the referee for his/her valuable comments that have led to an improvement of the manuscript.

\bibliographystyle{plain}

\end{document}